\documentclass[a4paper,12pt,reqno]{amsart}
\usepackage[utf8]{inputenc}

\usepackage{amsthm, amsmath, amssymb}
\usepackage{graphicx}
\usepackage{sidecap}
\usepackage{subfig}
\usepackage{comment}
\usepackage{color}
\usepackage{hyperref}
\hypersetup{
colorlinks=true,
linkcolor=red,
urlcolor=blue,
citecolor=red
}
\newtheorem{theorem}{Theorem}

\newtheorem{lemma}{Lemma}

\date{}

\author{Avetik Arakelyan}
\address{Institute of Mathematics, NAS of Armenia and Yerevan State University, Yerevan, Armenia}
\email{avetik.arakelyan@ysu.am}

\title[Convergence  of a Numerical Algorithm  for Two Population Densities] {Convergence Analysis of a Numerical Algorithm for Spatial Segregation of Two Population Densities}


\begin{document}

\begin{abstract}
Over the last few decades, the numerical approximation of spatial segregation equations for reaction-diffusion systems with $m$ population densities has gained significant interest. These problems are governed by a minimization problem subject to a closed but non-convex set. In the present work, we deal with the numerical approximation of equations of stationary states for a certain class of the spatial segregation of reaction-diffusion systems with two population densities having disjoint support. We prove the convergence of the numerical algorithm for two competing populations with non-negative internal dynamics $f_i(x)\geq 0$. At the end of the paper, we present computational tests.
\end{abstract}

\maketitle

\section{Introduction}

\subsection{The statement of the problem}
In recent decades, spatial segregation in reaction-diffusion systems has been intensively studied. The existence and uniqueness of spatially inhomogeneous solutions for Lotka-Volterra competition models with two or more competing densities have been investigated extensively 
\cite{arakelyan2017uniqueness,MR2146353,MR2151234,MR2283921,MR2300320, MR1900331, MR2417905}.
The aim of this paper is to study numerical solutions for a specific class of spatial segregation in reaction-diffusion systems with two population densities. This problem is closely related to systems with an arbitrary number of competing densities, which are governed by a minimization problem over a closed but non-convex set.

Let $\Omega \subset \mathbb{R}^n$ ($n\geq 2$) be a bounded, connected domain with a smooth boundary, and let $m \ge 2$ be a fixed integer. We consider the steady-states of $m$ competing species coexisting in the same area $\Omega$. Let $u_{i}(x)$ denote the population density of the $i^\textrm{th}$ component with the internal dynamic prescribed by $f_{i}(x)$. Here, we assume that $f_i$ is uniformly continuous and $f_{i}(x) \geq 0.$

An $m$-tuple $U=(u_1,\dots,u_m)\in (H^{1}(\Omega))^{m}$ is called a \emph{segregated state} if
\[
u_{i}(x) \cdot  u_{j}(x)=0, \quad \text{a.e. in } \Omega \text{ for } i\neq j.
\]
The problem amounts to minimizing the energy functional
\begin{equation}\label{1}
    E(u_1, \dots, u_m)=\int_{\Omega} \sum_{i=1}^{m} \left( \frac{1}{2}| \nabla u_{i}|^{2}+f_iu_i \right) dx,
\end{equation}
over the set
$$
S = \left\{ (u_1,\dots,u_{m})\in (H^{1}(\Omega))^{m} : u_{i}\geq 0, \, u_{i} \cdot u_{j}=0, \, u_{i}=\phi_{i} \text{ on } \partial \Omega \right\},
$$
where $\phi_{i} \in H^{\frac{1}{2}}(\partial \Omega),$ $\phi_{i} \cdot \phi_{j}=0$ for $i\neq j,$ and $\phi_{i}\geq 0$ on the boundary $\partial \Omega.$

Throughout this paper, we will focus on the case $m=2.$ The minimization problem is thus reduced to:
\begin{equation}\label{101}
    \text{Minimize } E(u_1, u_2)=\int_{\Omega} \sum_{i=1}^{2} \left( \frac{1}{2}| \nabla u_{i}|^{2}+f_iu_i\right) dx,
\end{equation}
over the set
$$
S = \left\{ (u_1,u_{2})\in (H^{1}(\Omega))^{2} : u_{i}\geq 0, \, u_{1} \cdot u_{2}=0, \, u_{i}=\phi_{i} \text{ on } \partial \Omega \right\}.
$$
Here, $\phi_{i} \in H^{\frac{1}{2}}(\partial \Omega)$ satisfies the properties $\phi_{1} \cdot \phi_{2}=0$ and $\phi_{i}\geq 0$ on the boundary $\partial \Omega.$

Due to the non-convexity of the set $S,$  the convergence analysis of the numerical scheme becomes more difficult.

\subsection{Two-phase membrane (obstacle) problem}
In this section, we briefly introduce the two-phase membrane problem and demonstrate its connection to the segregation problem for two competing densities (detailed discussions can be found in \cite{Mywork}). This connection plays a key role in proving the convergence of the proposed algorithm.

Let $f_i:\Omega \rightarrow \mathbb{R}$ for $i=1,2$ be non-negative Lipschitz continuous functions, where $\Omega$ is a bounded open subset of $\mathbb{R}^n$ with a smooth boundary. Let
\[
K=\{ v \in W^{1,2}(\Omega): v-g \in W^{1,2}_{0}(\Omega) \},
\]
where $g$ changes sign on the boundary. Consider the functional
\begin{equation}\label{13}
    I(v)=\int_{\Omega}\left(\frac{1}{2}|\nabla v|^2+f_{1} \max(v,0)-f_{2}\min(v,0)\right)dx,
\end{equation}
which is convex and weakly lower semi-continuous, and hence attains its unique infimum at some point $u \in K$.

The Euler-Lagrange equation corresponding to the minimizer $u$, commonly referred to as the two-phase obstacle problem, is given in \cite{MR1620644} and reads as follows
\begin{equation}\label{14}
\begin{cases}
    \Delta u = f_{1} \chi_{\{u >0 \}}-f_{2} \chi_{\{u <0\}}  &  \text{in } \Omega, \\
    u = g      & \text{on } \partial \Omega,
\end{cases}
\end{equation}
where $\Gamma(u) =\partial \{ x \in \Omega: u(x)>0 \} \cup \partial\{x\in \Omega: u(x) < 0\} \cap \Omega$ is called the \emph{free boundary.}

We define the positive and negative parts as
\[
v\vee 0=\max(v,0),\quad v\wedge 0=\min(v,0). 
\]

In the functional (\ref{13}), we set
\begin{align*}
    u_{1} = v\vee 0, \quad   u_{2} = -v\wedge 0, \\
    g_{1} = g\vee 0, \quad   g_{2} = -g\wedge 0.
\end{align*}

Then, the functional $I(v)$ in \eqref{13} can be rewritten as
\begin{equation}\label{functional}
    I(u_1,\,u_2)=\int_{\Omega}\left(\frac{1}{2}|\nabla u_{1}|^2+\frac{1}{2}|\nabla u_{2}|^{2}+ f_{1}u_1+ f_{2}u_{2} \right)dx,
\end{equation}
where the minimization is evaluated over the set
$$
S = \left\{ (u_1,\,u_{2})\in (H^{1}(\Omega))^{2} : u_{1} \cdot u_{2}=0,\, u_{i}\geq 0,\, u_{i}=g_{i} \text{ on } \partial \Omega, \, i=1,2 \right\}.
$$

If we set $u=u_1-u_2$ in system \eqref{14}, we obtain:
\begin{equation}\label{15}
\begin{cases}
    \Delta (u_1-u_2)= f_{1} \chi_{\{u_1-u_2 >0 \}}-f_{2} \chi_{\{u_1-u_2 <0\}}  &  \text{in } \Omega, \\
    u_1-u_2=g_1-g_2      & \text{on } \partial \Omega.
\end{cases}
\end{equation}
Thus, the solution to our minimization problem \eqref{101} satisfies the two-phase obstacle problem in the distributional sense as written in \eqref{15}. In the case of three or more competing densities, this property is fulfilled only locally  \cite[Proposition 6.3]{MR2151234}.

For more information regarding the two-phase membrane (obstacle) problem, the interested reader is referred to the book \cite{PSU2012} and the references therein. For numerical analysis, we refer to \cite{MR2961456,WithRafMic,Farid}, and for the parabolic version of this problem, we direct the reader to \cite{arakelyan2015finite}. Recently, optimal control configurations governed by the two-phase membrane problem have also been explored \cite{bozorgnia2025optimal}. In these models, the non-differentiable control-to-state mapping requires regularization of the state equation to establish the existence, uniqueness, and characterization of optimal control properties.

\subsection{Known results}
Over the last few decades, there has been considerable interest in studying the numerical approximation of reaction-diffusion equations, particularly those arising in population ecology when highly competitive interactions between different species occur.

We refer the reader to \cite{MR2151234,MR2363653, MR2079274, MR1303035, MR1687440, MR1900331}  and in particular to \cite{MR1687440} for models involving Dirichlet boundary data. A complete analysis of the stationary case has been studied  in \cite{MR2151234}. Also  numerical simulation   for the spatial segregation limit of two diffusive Lotka-Volterra models in presence of strong competition and inhomogeneous Dirichlet boundary conditions is provided in \cite{MR2459673}. The authors in \cite{MR2459673}  solve the problem for small  $\varepsilon$ and then let $\varepsilon  \longrightarrow 0.$

 In the work \cite{MR2563520} Bozorgnia proposed two numerical algorithms for the problem \eqref{1} with lack of internal dynamics ($f_i=0$). The finite element approximation is based on the local properties of the solution. In this case the author was able to provide the convergence of the method. Unfortunately, this nice idea cannot be generalized for the case with non-negative internal dynamics. The second approach is a finite difference method, but lack of analysis of the scheme. This finite difference method has been generalized in \cite{Mywork} for the case of non-negative $f_i$. In \cite{Mywork} the authors give a numerical consistent variational system with strong interaction, and provide disjointness condition of populations during the iteration of the scheme.

 The further advancements have been made in developing and analyzing finite difference approaches for these problems. In \cite{arakelyan2016numerical} the authors proposed a generalized numerical approach for the spatial segregation of reaction-diffusion systems with $m$ population densities. They prove the existence and uniqueness of appropriate finite difference scheme. By using the developed tools in \cite{arakelyan2016numerical}, the current author was able   to establish the convergence of the finite difference scheme for a general class of these spatial segregation systems \cite{arakelyan2018convergence}.

 The present work deals with the analysis of the convergence of the algorithm  for two competing populations with non-negative  internal dynamics $f_i.$

\subsection{Notations}
We will define the notations for the one-dimensional and two-dimensional cases in parallel.

For the sake of simplicity, we will assume that $\Omega=(-1,1)$ in one-dimensional case and $\Omega=(-1,1)\times(-1,1)$ in two-dimensional case in the rest of the paper, keeping in mind that the method works also for more complicated domains.

Let $N\in\mathbb{N}$ be a positive integer, $h=2/N$ and
$$
x_i=-1+ ih,\, y_i=-1+ih,\quad i=0,1,...,N.
$$

We use the notation $u_i$ and $u_{i,j}$ (or simply $u_\alpha$, where $\alpha$ is one- or two-dimensional index) for finite-difference scheme approximation to $u(x_i)$ and $u(x_i,y_j)$,
\[f_1(i)=f_1(x_i),\;\; f_2(i)=f_2(x_i),\]
\[f_1(i,j)=f_1(x_i,y_j),\;\;f_2(i,j)=f_2(x_i,y_j),\]
\[g_i=\phi_1(i)-\phi_2(i)=\phi_1(x_i)-\phi_2(x_i)\]
and
\[g_{i,j}=\phi_1(i,j)-\phi_2(i,j)=\phi_1(x_i,y_j)-\phi_2(x_i,y_j),\]
in one- and two-dimensional cases, respectively, assuming that the function $\phi_1-\phi_2$  is extended to be zero everywhere outside the boundary $\partial\Omega$.

In this paper we will use also notations $u=(u_{\alpha})$, $g=(g_{\alpha})$ (not to be confused with functions $u, g$).

Denote
$$
 \mathcal N=\{i:\ 0\leq i\leq N\}\quad\mbox{or}\quad \mathcal N=\{(i,j):\ 0\leq i,j\leq N\},
$$
$$
  \mathcal N^o=\{i:\ 1\leq i\leq N-1\} \quad\mbox{or}\quad  \mathcal N^o=\{(i,j):\ 1\leq i,j\leq N-1\},
$$
in one- and two- dimensional cases, respectively, and
$$
  \partial \mathcal N=\mathcal N \setminus \mathcal N^o.
$$

In one-dimensional case we consider the following approximation for Laplace operator: for any $i\in \mathcal N^o$,
\[
  \Delta_h u_{i}\equiv L_h u_{i}=\frac{u_{i-1}-2u_{i}+u_{i+1}}{h^2},
\]
and for two-dimensional case we introduce the following 5-point stencil approximation for Laplacian:
\[
 \Delta_h u_{i,j}\equiv L_h u_{i,j}=\frac{u_{i-1,j}+u_{i+1,j}-4u_{i,j}+u_{i,j-1}+u_{i,j+1}}{h^2}
\]
for any $(i,j)\in \mathcal N^o$.


\section{Numerical algorithm and its properties}
In this section we consider an algorithm, which is basically the generalization of the numerical algorithm developed in \cite{Mywork}, for the case $f_i(x,s)=f_i(x).$ The iterative method for two densities, say $u_{1}$ and $u_{2},$ will be as follows:

$\bullet$ \textbf{Initialization:}
\begin{equation*}
u_{1}^{0}(x_{i},y_{j})=
\left \{
\begin{array}{ll}
  0    &    (x_{i},y_{j}) \in \mathcal N^o,   \\
  \phi_{1}(x_{i},y_{j})   &   (x_{i},y_{j}) \in \partial \mathcal N.
\end{array}
\right.
\end{equation*}
\begin{equation*}
u_{2}^{0}(x_{i},y_{j})=
\left \{
\begin{array}{ll}
  0                     &   (x_{i},y_{j}) \in \mathcal N^o,     \\
  \phi_{2}(x_{i},y_{j})  &     (x_{i},y_{j}) \in \partial\mathcal N.
\end{array}
\right.
\end{equation*}

$\bullet$ \textbf{Step $k+1$, $k\geq 0$:}

We iterate over all interior points by setting:
\begin{equation}\label{algorithm}
\begin{cases}
u_1^{k+1} (x_{i},y_{j}) = \max\left(\frac{- f_1(x_i,y_j) h^2}{4}+\overline{u}_1^{k}(x_{i},y_{j})-\overline{u}_2^{k}(x_{i},y_{j}), 0\right),\\[8pt]
u_2^{k+1} (x_{i},y_{j}) =\max\left( \frac{- f_{2}(x_i,y_j) h^2}{4}+\overline{u}_2^{k}(x_{i},y_{j})-\overline{u}_1^{k}(x_{i},y_{j}), 0\right).
\end{cases}
\end{equation}

Here for a given uniform mesh on $\Omega\subset \mathbb{R}^2,$ we define $\overline{u}_l^{k}(x_{i},y_{j})$ for $l=1,2,$ to be the average of $u_l$ for all neighbor points of $(x_{i},y_{j})$ using the forward Gauss-Seidel update
\[
\overline{u}_l^{k}(x_{i},y_{j})=\frac{1}{4}\left[u_l^{k+1} (x_{i-1},y_{j})+u_l^{k} (x_{i+1},y_{j})+u_l^{k+1}(x_{i},y_{j-1})+u_l^{k} (x_{i},y_{j+1})\right].
\]
Next, we discuss some properties of the given algorithm.
\begin{lemma}\label{lemma_1}
The numerical algorithm \eqref{algorithm} satisfies
\[
u_{1}^{k}(x_i,y_j)\cdot u_{2}^{k}(x_i,y_j)=0,
\]
for all $k\in\mathbb{N}.$
\end{lemma}
\begin{proof}
    The proof of the disjoint property of the densities in the general case can be found in \cite{Mywork}. Hence, we omit the proof.
\end{proof}

\begin{lemma}\label{lemma_2}
The numerical algorithm \eqref{algorithm} is stable and consistent.
\end{lemma}

\begin{proof}
Here we will prove the stability of the method. For the proof of consistency, we refer the reader to \cite{Mywork}.

Due to the non-negative $f_i\geq 0,$ we can write the following inequalities:
\[
u_1^{k+1} (x_{i},y_{j}) = \max\left(\frac{- f_1(x_i,y_j) h^2}{4}+\overline{u}_{1}^{k}(x_{i},y_{j})-\overline{u}_{2}^{k}(x_{i},y_{j}),\quad 0\right)\leq \overline{u}_{1}^{k}(x_{i},y_{j}),
\]
and
\[
u_2^{k+1} (x_{i},y_{j}) = \max\left(\frac{- f_2(x_i,y_j) h^2}{4}+\overline{u}_{2}^{k}(x_{i},y_{j})-\overline{u}_{1}^{k}(x_{i},y_{j}),\quad 0\right)\leq \overline{u}_{2}^{k}(x_{i},y_{j}).
\]
Therefore
\[
u_1^{k+1} (x_{i},y_{j})\leq \overline{u}_{1}^{k}(x_{i},y_{j}),\;\; \text{and}\;\; u_2^{k+1} (x_{i},y_{j})\leq \overline{u}_{2}^{k}(x_{i},y_{j}),
\]
respectively.

We apply the discrete maximum principle via induction. Let
$$
M_l = \max_{(x_i,y_j)\in\partial \mathcal{N}} \phi_{l}(x_{i},y_{j}).
$$
For $k=0$, the initialization provides
$$
0 \leq u_l^0(x_i,y_j) \leq M_l\;\;\text{for all}\;\;(x_i,y_j) \in \mathcal{N}.
$$

Assume $u_l^k(x_i,y_j) \leq M_l$ globally. Because the Gauss-Seidel algorithm updates nodes sequentially, the newly computed neighbor values $u_l^{k+1}$ used in the average already satisfy the bound $M_l$, i.e. $u_l^{k+1}(x_{i-1},y_{j})\le M_l$ and  $u_l^{k+1}(x_{i},y_{j-1})\le M_l$. Thus,
\[
\overline{u}_l^{k}(x_{i},y_{j}) = \frac{1}{4}\left[ u_l^{k+1}(x_{i-1},y_{j}) + u_l^{k}(x_{i+1},y_{j}) + u_l^{k+1}(x_{i},y_{j-1}) + u_l^{k}(x_{i},y_{j+1})\right] \leq M_l.
\]
It implies that
\[
u_l^{k+1}(x_{i},y_{j}) \leq \overline{u}_l^{k}(x_{i},y_{j}) \leq M_l.
\]
By mathematical induction, we conclude that for all points $(x_i, y_j) \in \mathcal{N}$:
\[
0 \leq u_1^{k+1}(x_{i},y_{j}) \leq \max_{(x_i,y_j)\in\partial \mathcal{N}}\phi_{1}(x_{i},y_{j}),
\]
\[
0 \leq u_2^{k+1}(x_{i},y_{j}) \leq \max_{(x_i,y_j)\in\partial \mathcal{N}}\phi_{2}(x_{i},y_{j}).
\]
Hence, the sequences $u_1^{k}$ and $u_2^{k}$ are uniformly bounded for every $k\in\mathbb{N}$. This completes the proof of stability.
\end{proof}
\section{Algorithm convergence}

In this section we show convergence of the numerical algorithm for both one- and two-dimensional cases. 

\subsection{Algorithm for one-dimensional case}

For the sake of simplicity, we consider here only the one-dimensional case. Let $u=(u_0,u_1,\dots,u_N)$ be the solution of \eqref{algorithm} in one-dimensional case. In this case the algorithm reads:

$\bullet$ \textbf{Initialization:}
 \begin{equation*}
 u_{1}^{0}(x_{i})=
\left \{
\begin{array}{ll}
  0    &    x_{i} \in \mathcal{N}^o,   \\
  \phi_{1}(x_{i})   &   x_{i} \in \partial \mathcal{N}.
\end{array}
\right.
\end{equation*}
\begin{equation*}
u_{2}^{0}(x_{i})=
\left \{
\begin{array}{ll}
  0                     &   x_{i} \in \mathcal{N}^o,     \\
  \phi_{2}(x_{i})  &     x_{i} \in \partial \mathcal{N}.
\end{array}
\right.
\end{equation*}

$\bullet$ \textbf{{Step $k+1$, $k\geq 0:$}}

We iterate over all interior points by setting
\begin{equation}\label{algorithm1D}
\begin{cases}
u_1^{k+1} (x_{i}) = \max\left(\frac{- f_1(x_i) h^2}{2}+\overline{u}_{1}^{k}(x_{i})-\overline{u}_{2}^{k}(x_{i}), 0\right),\\[8pt]
u_2^{k+1} (x_{i}) =\max\left( \frac{- f_{2}(x_i) h^2}{2}+\overline{u}_{2}^{k}(x_{i})-\overline{u}_{1}^{k}(x_{i}), 0\right).
\end{cases}
\end{equation}
Here for a given uniform mesh on $\Omega\subset \mathbb{R},$ we define $\overline{u}_l^{k}(x_{i})$ for $l=1,2,$ to be the average of $u_l$ for all neighbor points of $x_{i}$ using the Gauss-Seidel update
\[
\overline{u}_l^{k}(x_{i})=\frac{1}{2}[u_l^{k+1} (x_{i-1})+u_l^{k} (x_{i+1})].
\]

Define $g(x)=\phi_1(x)-\phi_2(x)$ extended by zero on the interior $\mathcal{N}^o$, and let $w_i\equiv {w}(x_{i})={u_1}(x_{i})-{u_2}(x_{i})$ be the physical density difference. We consider the decomposition $w = v + g$, where $v$ captures the interior variation. 

We consider the following discrete functional evaluated over $v$:
\[
J_h(v)=-\frac{1}{2} \Big(L_h v,v\Big)+\Big(f_1, v \vee 0\Big)-\Big(f_2,v\wedge 0\Big)-\Big(L_h g,v\Big),
\]
defined on the finite dimensional space
\[
\mathcal K=\{v\in\mathcal H: \ v_{\alpha}=0, \ \alpha\in\partial \mathcal {N}\},\quad {\rm where}\quad \mathcal H=\{v=(v_\alpha): v_\alpha\in \mathbb {R}, \ \alpha\in \mathcal {N}\}.
\]
Here $v\vee 0=\max(v,0)$, $v\wedge 0=\min(v,0)$ and for two sequences $(a_\alpha)$ and $(b_\alpha)$, $\alpha\in \mathcal N$, the inner product $(\cdot, \cdot)$ is defined by
\[
(a,b)=\sum_{\alpha\in \mathcal N}a_\alpha\cdot b_\alpha.
\]
In particular, since $w=v+g\in\mathcal{H}$ and $v\in\mathcal{K},$ we recover the exact boundary conditions $w_0=g_0$ and $w_N=g_N$. Furthermore, because $g=0$ on $\mathcal{N}^o$, it follows that $v_i = w_i = u_1(x_i) - u_2(x_i)$ for all interior nodes $i \in \mathcal{N}^o$.

We will use the notation $\tilde v=(v_1,v_2,\dots,v_{N-1})$. This is the unknown part in $v$ that needs to be calculated. 

In the next section we will prove the convergence of $\tilde v^{k}=(v_1^{k},v_2^{k},\dots,v_{N-1}^{k}),$ and then the disjointness condition for competing densities will lead to the convergence of $u_1^{k}$ and $u_2^{k}$ separately.

\subsection{Convergence result for the one-dimensional case}

\begin{theorem}\label{2P-thm-PGS-convergence}
The sequence $\tilde v^k$ converges.
\end{theorem}

\begin{proof}
Denote the intermediate update states as:
\[
\tilde v^{k,i}=\left(\tilde v^k_1,\tilde v^k_2,\dots,\tilde v^k_i, \tilde v^{k-1}_{i+1},\dots,\tilde v^{k-1}_{N-1}\right),\quad i=1,\dots,N-1,\quad k\in\mathbb N.
\]
We embed $\tilde{v}^{k,i}$ into the full space $\mathcal{K}$ to form $v^{k,i}$ by fixing the boundary nodes to $0$:
\[
v^{k,i}=\left(0,\tilde v^k_1,\tilde v^k_2,\dots,\tilde v^k_i, \tilde v^{k-1}_{i+1},\dots,\tilde v^{k-1}_{N-1},0\right)\in\mathcal K,\quad i=1,\dots,N-1,\quad k\in\mathbb N.
\]
We evaluate the discrete functional $\mathcal J_p= J_h\left(v^{k,i}\right)$ for the global step index $p=(N-1)(k-1)+i$ with $i=1,\dots,N-1.$

The main idea is to prove that $\mathcal J_p$ decreases.

First let $p\notin\{q(N-1): q\in\mathbb N\}$, meaning $i\ne N-1$. The transition from $\mathcal J_p$ to $\mathcal J_{p+1}$ represents the update of a single node corresponding to the sequential index $i+1$. Let $v^{old} = \tilde v^{k-1}_{i+1}$ and $v^{new} = \tilde v^{k}_{i+1}$.

Because the diagonal entry of the 1D negative discrete Laplacian $L_h$ is $\frac{2}{h^2}$, the algebraic expansion of the energy difference yields:
\begin{multline*}   
\mathcal J_p-\mathcal J_{p+1}= J_h\left(v^{k,i}\right)- J_h\left(v^{k,i+1}\right)=-\frac{1}{2}\left(L_h \left(v^{k,i}-v^{k,i+1}\right),v^{k,i}-v^{k,i+1}\right)-\\-
\left(L_h v^{k,i+1},v^{k,i}-v^{k,i+1}\right)+\left(f_1,v^{k,i}\vee 0-v^{k,i+1}\vee 0\right)-\left(f_2,v^{k,i}\wedge 0-v^{k,i+1}\wedge 0\right)-\\-\left(L_h g,v^{k,i}-v^{k,i+1}\right)
=\frac{1}{h^2}\left(v^{old}-v^{new}\right)^2
-(L_h v^{k,i+1})_{i+1}\cdot \left(v^{old}-v^{new}\right)+\\+f_1(i+1)\cdot\left[v^{old}\vee 0-v^{new}\vee 0\right]
-f_2(i+1)\cdot\left[v^{old}\wedge 0-v^{new}\wedge 0\right]-(L_h g)_{i+1}\cdot \left(v^{old}-v^{new}\right).
\end{multline*}

We continue by considering three cases for the newly updated value $v^{new}$:
\\

\paragraph{\textbf{Case 1: $v^{new} > 0$}} Due to $u_1^k(x_{i+1})\cdot u_2^k(x_{i+1})=0,$ we have $ v^{new}=u_1^k(x_{i+1})>0,$ and $u_2^k(x_{i+1})=0.$ It follows from the Gauss-Seidel update \eqref{algorithm1D} that the discrete Laplacian exactly balances the right-hand side and boundary terms:
\[
-(L_h v^{k,i+1})_{i+1} - (L_h g)_{i+1} = -f_1(i+1).
\]
Substituting this equality into the expansion, the linear terms cancel perfectly, yielding:
\begin{align*}
  \mathcal J_p-\mathcal J_{p+1}& = \frac{1}{h^2}\left(v^{old}-v^{new}\right)^2 - f_1(i+1)\cdot \left(v^{old}-v^{new}\right) +\\&+ f_1(i+1)\cdot\left(v^{old}\vee 0-v^{new}\right)- f_2(i+1)\cdot (v^{old}\wedge 0).  
\end{align*}
Simplifying this expression gives
\[
\mathcal J_p-\mathcal J_{p+1}=\frac{1}{h^2}\left(v^{old}-v^{new}\right)^2 - (f_1(i+1)+f_2(i+1))\cdot \left(v^{old}\wedge 0\right) \ge \frac{1}{h^2}\left(v^{old}-v^{new}\right)^2.
\]
Hence, in this case we have
\begin{equation}\label{2p-eq-thm3.1-1}
\mathcal J_p-\mathcal J_{p+1}\ge\frac{1}{h^2}\left(\tilde v^{k-1}_{i+1}-\tilde v^{k}_{i+1}\right)^2.
\end{equation}
\\

\paragraph{\textbf{Case 2: $v^{new} < 0$}} In this case again due to $u_1^k(x_{i+1})\cdot u_2^k(x_{i+1})=0,$ we have $v^{new}=-u_2^k(x_{i+1})<0,$ and $u_1^k(x_{i+1})=0.$ Analogously to the previous case, the discrete Laplacian perfectly balances $f_2(i+1)$, and we can prove that \eqref{2p-eq-thm3.1-1} holds also in this case.
\\

\paragraph{\textbf{Case 3: $v^{new}=0$}} Since $v^{new}=u_1^k(x_{i+1})-u_2^k(x_{i+1})=0,$ then disjointness implies $u_1^k(x_{i+1})=u_2^k(x_{i+1})=0.$ Thus, according to \eqref{algorithm1D}, the discrete Laplacian is bounded by the non-negative functions $f_1(i+1)$ and $f_2(i+1)$:
\[
-f_1(i+1) \leq -(L_h v^{k,i+1})_{i+1} - (L_h g)_{i+1} \leq f_2(i+1).
\]
Substituting these bounds based on whether $v^{old}$ is positive or negative, we deduce that \eqref{2p-eq-thm3.1-1} holds in this case as well.

So far we have considered the case $p\notin\{q(N-1): q\in\mathbb N\}$. Now assume that $p\in\{q(N-1): q\in\mathbb N\}$. This represents the transition between sweeps, updating the very first node $i=0$ (evaluating for physical node $1$). In this case, identical logic yields:
\begin{equation}\label{2p-eq-thm3.1-2}
\mathcal J_p-\mathcal J_{p+1}\ge \frac{1}{h^2}\left(\tilde v^{k+1}_{1}-\tilde v^{k}_{1}\right)^2.
\end{equation}

Summarizing, we deduce that $\mathcal J_p$ decreases, and, since it is bounded from below, the sequence $\mathcal J_p$ converges. By summing the inequalities \eqref{2p-eq-thm3.1-1} and \eqref{2p-eq-thm3.1-2} over all steps, we obtain that $\sum (\tilde v^k_i - \tilde v^{k-1}_i)^2 < \infty$, which implies $|\tilde v^k_i - \tilde v^{k-1}_i| \to 0$ as $k \to \infty$. 
Furthermore, because the discrete functional $J_h(v)$ is coercive, the sequence of iterates $\tilde{v}^k$ is uniformly bounded, ensuring the existence of at least one limit point $\tilde{v}^*$. Passing to the limit in the discrete Projected Gauss-Seidel updates reveals that any such limit point $\tilde{v}^*$ must be a fixed point of the iteration, meaning it perfectly satisfies the coordinate-wise optimality conditions of the minimization problem. Due to the strict convexity of $J_h(v)$, this stationary point is the unique global minimizer. Consequently, every convergent subsequence must converge to this same unique minimizer $\tilde{v}^*$, meaning the entire bounded sequence $\tilde{v}^k$ converges. Since it converges, $\tilde{v}^k_i$ is a Cauchy sequence for every fixed $i=1,\dots,N-1$.
\end{proof}

Thus we have proved that the sequence $\tilde v^{k}_{i}=u_1^k(x_{i})-u_2^k(x_{i})$ is converging to $\tilde v$ for every fixed $i=1,\dots,N-1.$ Observe that due to the disjointness property, we have
\begin{equation}\label{disjoint}
u_1^k(x_{i})=\max(u_1^k(x_{i})-u_2^k(x_{i}),0)\;\;\mbox{and}\;\; u_2^k(x_{i})=\max(u_2^k(x_{i})-u_1^k(x_{i}),0).
\end{equation}
Now, in view of \eqref{disjoint} we will obtain the convergence of $u_1^k(x_{i})$ and $u_2^k(x_{i})$ separately.

\subsection{Convergence result for the two-dimensional case}

To formalize the proof for the two-dimensional algorithm, we map the 2D grid indices $(i,j) \in \mathcal{N}^o$ to a 1D sequence using a standard row-wise ordering. Let $N_{int} = (N-1)^2$ be the total number of interior nodes, and let the bijection $m(i,j) \in \{1, \dots, N_{int}\}$ denote the sequential index of the node $(i,j)$. We define the $N_{int}$-dimensional vector $\tilde{v}^k = (v_1^k, v_2^k, \dots, v_{N_{int}}^k)$, where $v_{m(i,j)}^k = u_{1,i,j}^k - u_{2,i,j}^k$.

\begin{theorem}\label{2P-thm-PGS-convergence-2D}
The sequence $\tilde v^k$ converges.
\end{theorem}

\begin{proof}
Denote the intermediate update states for the interior nodes as:
\[
\tilde v^{k,m}=\left(\tilde v^k_1, \dots, \tilde v^k_m, \tilde v^{k-1}_{m+1}, \dots, \tilde v^{k-1}_{N_{int}}\right), \quad m=1,\dots,N_{int},\quad k\in\mathbb N.
\]
We embed $\tilde{v}^{k,m}$ into the full space $\mathcal{K}$ to form $v^{k,m}$ by fixing the boundary nodes to $0$. The physical boundary data $g_{i,j}$ is natively handled by the $-(L_h g, v)$ term in the functional. We evaluate the discrete functional $\mathcal J_p = J_h\left(v^{k,m}\right)$ for the global step index $p = (k-1)N_{int} + m$.

The main idea is to prove that $\mathcal J_p$ decreases.

First, let $p \notin \{q N_{int} : q \in \mathbb{N}\}$, meaning $m \ne N_{int}$. The transition from $\mathcal J_p$ to $\mathcal J_{p+1}$ represents the update of a single node $(i,j)$ corresponding to the sequential index $m+1$. Let $v^{old} = \tilde v^{k-1}_{m+1}$ and $v^{new} = \tilde v^{k}_{m+1}$.

Because the diagonal entry of the 2D negative discrete Laplacian $L_h$ is $\frac{4}{h^2}$, the algebraic expansion of the energy difference yields
\begin{multline*}
    \mathcal J_p-\mathcal J_{p+1}= J_h\left(v^{k,m}\right)- J_h\left(v^{k,m+1}\right)=-\frac{1}{2}\left(L_h \left(v^{k,m}-v^{k,m+1}\right),v^{k,m}-v^{k,m+1}\right)-\\-
\left(L_h v^{k,m+1},v^{k,m}-v^{k,m+1}\right)+\left(f_1,v^{k,m}\vee 0-v^{k,m+1}\vee 0\right)-\left(f_2,v^{k,m}\wedge 0-v^{k,m+1}\wedge 0\right)-\\-\left(L_h g,v^{k,m}-v^{k,m+1}\right)
=\frac{2}{h^2}\left(v^{old} - v^{new}\right)^2
- (L_h v^{k,m+1})_{i,j} \cdot \left(v^{old}-v^{new}\right) +\\+ f_1(i,j)\cdot\left[v^{old}\vee 0-v^{new}\vee 0\right]
-f_2(i,j)\cdot\left[v^{old}\wedge 0-v^{new}\wedge 0\right]-(L_h g)_{i,j}\cdot \left(v^{old}-v^{new}\right).
\end{multline*}

We again continue by considering three cases for the newly updated value $v^{new}$:
\\

\paragraph{\textbf{Case 1: $v^{new} > 0$}} Due to $u_{1,i,j}^k \cdot u_{2,i,j}^k = 0,$ we have $v^{new}=u_{1,i,j}^k >0,$ and $u_{2,i,j}^k =0.$ It follows from the Gauss-Seidel update \eqref{algorithm} that the discrete Laplacian exactly balances the right-hand side and boundary terms
\[
-(L_h v^{k,m+1})_{i,j} - (L_h g)_{i,j} = -f_1(i,j).
\]
Substituting this equality into the expansion, the linear terms cancel perfectly, yielding
\begin{align*}
  \mathcal J_p-\mathcal J_{p+1}& = \frac{2}{h^2}\left(v^{old}-v^{new}\right)^2 - f_1(i,j)\cdot \left(v^{old}-v^{new}\right) +\\&+ f_1(i,j)\cdot\left(v^{old}\vee 0-v^{new}\right)- f_2(i,j)\cdot (v^{old}\wedge 0).  
\end{align*}
Simplifying this expression gives
\[
\mathcal J_p-\mathcal J_{p+1}=\frac{2}{h^2}\left(v^{old}-v^{new}\right)^2 - (f_1(i,j)+f_2(i,j))\cdot \left(v^{old}\wedge 0\right) \ge \frac{2}{h^2}\left(v^{old}-v^{new}\right)^2.
\]
Hence, in this case we have
\begin{equation}\label{2p-eq-thm3.2-1}
\mathcal J_p-\mathcal J_{p+1}\ge\frac{2}{h^2}\left(\tilde v^{k-1}_{m+1}-\tilde v^{k}_{m+1}\right)^2.
\end{equation}
\\
\paragraph{\textbf{Case 2: $v^{new} < 0$}} In this case again due to $u_{1,i,j}^k \cdot u_{2,i,j}^k=0,$ we have $v^{new}=-u_{2,i,j}^k <0,$ and $u_{1,i,j}^k=0.$ Analogously to the previous case, the discrete Laplacian perfectly balances $f_2(i,j)$, and we can prove that \eqref{2p-eq-thm3.2-1} holds also in this case.
\\
\paragraph{\textbf{Case 3: $v^{new}=0$}} Since $v^{new}=u_{1,i,j}^k-u_{2,i,j}^k=0,$ then disjointness implies $u_{1,i,j}^k=u_{2,i,j}^k=0.$ Thus, according to \eqref{algorithm}, the discrete Laplacian is bounded by the non-negative functions $f_1(i,j)$ and $f_2(i,j)$:
\[
-f_1(i,j) \leq -(L_h v^{k,m+1})_{i,j} - (L_h g)_{i,j} \leq f_2(i,j).
\]
Substituting these bounds based on whether $v^{old}$ is positive or negative, we deduce that \eqref{2p-eq-thm3.2-1} holds in this case as well.

So far we have considered the case $p\notin\{q N_{int} : q\in\mathbb N\}$. Now assume that $p\in\{q N_{int} : q\in\mathbb N\}$. This represents the transition between sweeps, updating the very first node $m=1$. In this case, identical logic yields:
\begin{equation}\label{2p-eq-thm3.2-2}
\mathcal J_p-\mathcal J_{p+1}\ge \frac{2}{h^2}\left(\tilde v^{k+1}_{1}-\tilde v^{k}_{1}\right)^2.
\end{equation}

Summarizing, we deduce that $\mathcal J_p$ decreases, and, since it is bounded from below, the sequence $\mathcal J_p$ converges. By summing the inequalities \eqref{2p-eq-thm3.2-1} and \eqref{2p-eq-thm3.2-2} over all steps, we obtain that $\sum_{k=1}^\infty (\tilde v^k_{i,j} - \tilde v^{k-1}_{i,j})^2 < \infty$, which implies  $|\tilde v^k_{i,j} - \tilde v^{k-1}_{i,j}| \to 0$ as $k \to \infty$. By the same arguments as  in the proof of one-dimensional case we conclude that 
$\tilde{v}^k_{i,j}$ is a Cauchy sequence for every fixed spatial coordinate $(i,j)\in\mathcal N^o$.
\end{proof}

\section{Numerical examples}

In this section, we present numerical simulations for two competing densities with different internal dynamics $f_i$. We consider the following minimization problem:
\begin{equation}\label{last}
\text{Minimize} \int_{\Omega} \sum_{i=1}^{2} \left( \frac{1}{2}| \nabla u_{i}|^{2}+f_{i}u_{i} \right) dx,
\end{equation}
over the set
\[
S = \left\{ (u_1,u_{2})\in (H^{1}(\Omega))^{2} : u_{i}\geq 0, \, u_{1} \cdot u_{2}=0, \, u_{i}=\phi_{i} \text{ on } \partial \Omega \right\}.
\]

In all numerical experiments, the spatial domain is discretized using a uniform grid with step size $h = 1/N$. The projected Gauss-Seidel iterations are terminated when the maximum pointwise difference between successive sweeps falls below a prescribed tolerance, $\max |u_l^{k+1} - u_l^k| < 10^{-6}$. Note that, throughout all computations, the disjointness condition $u^k_{1}(x_i,y_j) \cdot u^k_{2}(x_i,y_j) = 0$ is strictly preserved at every grid point and every iteration, empirically validating the algebraic properties established in Lemma \ref{lemma_1}.

\begin{figure}[!htbp]
    \centering
    \includegraphics[width=\textwidth]{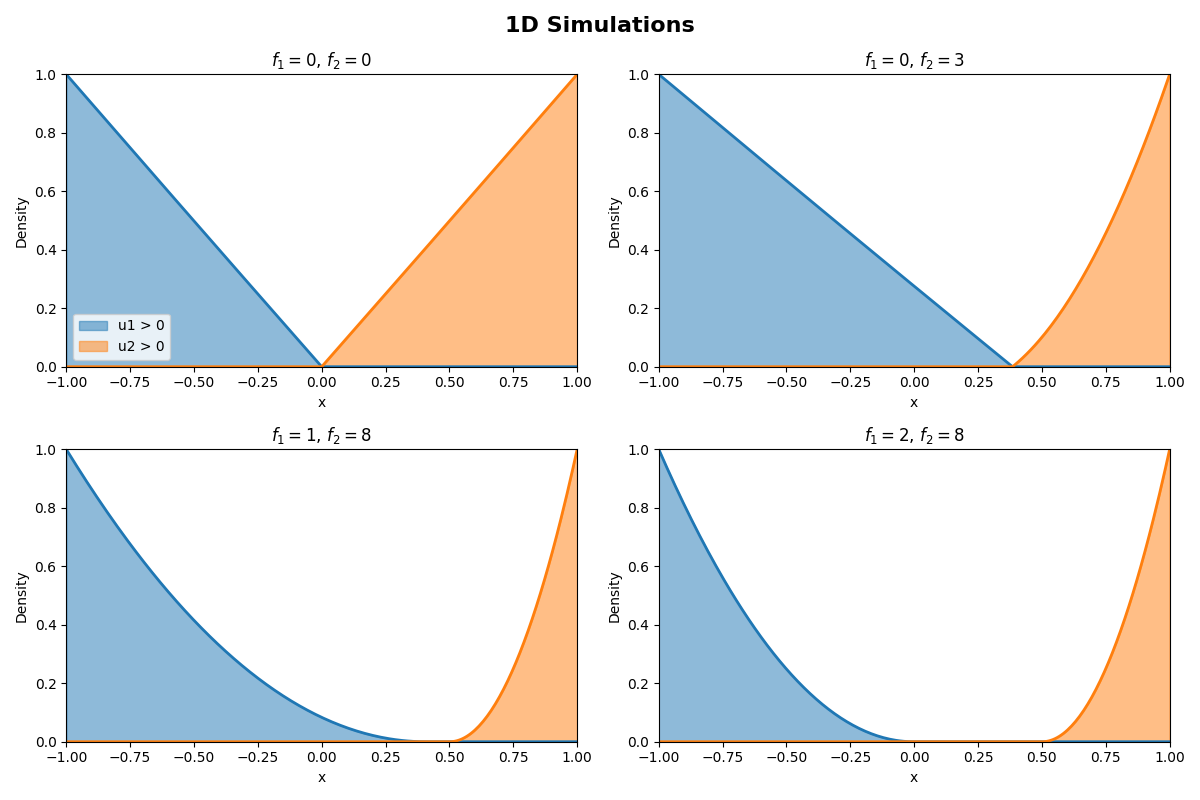}
 \caption{One-dimensional simulations with boundary conditions $\phi_1(-1)=\phi_2(1)=1$ and $\phi_2(-1)=\phi_1(1)=0$. The subfigures illustrate the effect of varying constant internal dynamics $f_i$, showing how larger penalties shrink the population support.}
  \label{1dfig}
\end{figure}

In Figure \ref{1dfig}, we consider the one-dimensional domain $\Omega=[-1,1]$, where the internal dynamics $f_1$ and $f_2$ are assumed to be constant. In the context of population dynamics, the parameter $f_i$ acts as a decay rate or environmental penalty. As expected, the species with the smaller internal dynamic $f_i$ faces less resistance, allowing it to occupy a larger portion of the domain. The emergence of the free boundaries is clearly visible, sharply defining the segregated regions.

Next, we present two-dimensional numerical experiments on the square domain $\Omega=[0,1]\times[0,1]$. The boundary conditions $\phi_1(x,y)$ and $\phi_2(x,y)$ are strictly disjoint and defined as follows:

\begin{equation*}
\phi_1(x,0) =
\begin{cases}
  0.5-2.5x  & 0 \leq x \leq 0.2,\\
  0         & 0.2 < x \leq 1,
\end{cases}
\qquad
\phi_1(x,1) =
\begin{cases}
  0.5-\frac{5}{8}x & 0 \leq x \leq 0.8,\\
  0                & 0.8 < x \leq 1,
\end{cases}
\end{equation*}
\[
\phi_1(0,y)=0.5, \qquad \phi_1(1,y)=0,
\]
and
\begin{equation*}
\phi_2(x,0) =
\begin{cases}
  0                         & 0 \leq x \leq 0.2,\\
  -\frac{1}{8}+\frac{5}{8}x & 0.2 < x \leq 1,
\end{cases}
\qquad
\phi_2(x,1) =
\begin{cases}
  0         & 0 \leq x \leq 0.8,\\
  -2+2.5x   & 0.8 < x \leq 1,
\end{cases}
\end{equation*}
\[
\phi_2(0,y)=0, \qquad \phi_2(1,y)=0.5.
\]

\begin{figure}[!htbp]
    \centering
    \includegraphics[width=\textwidth]{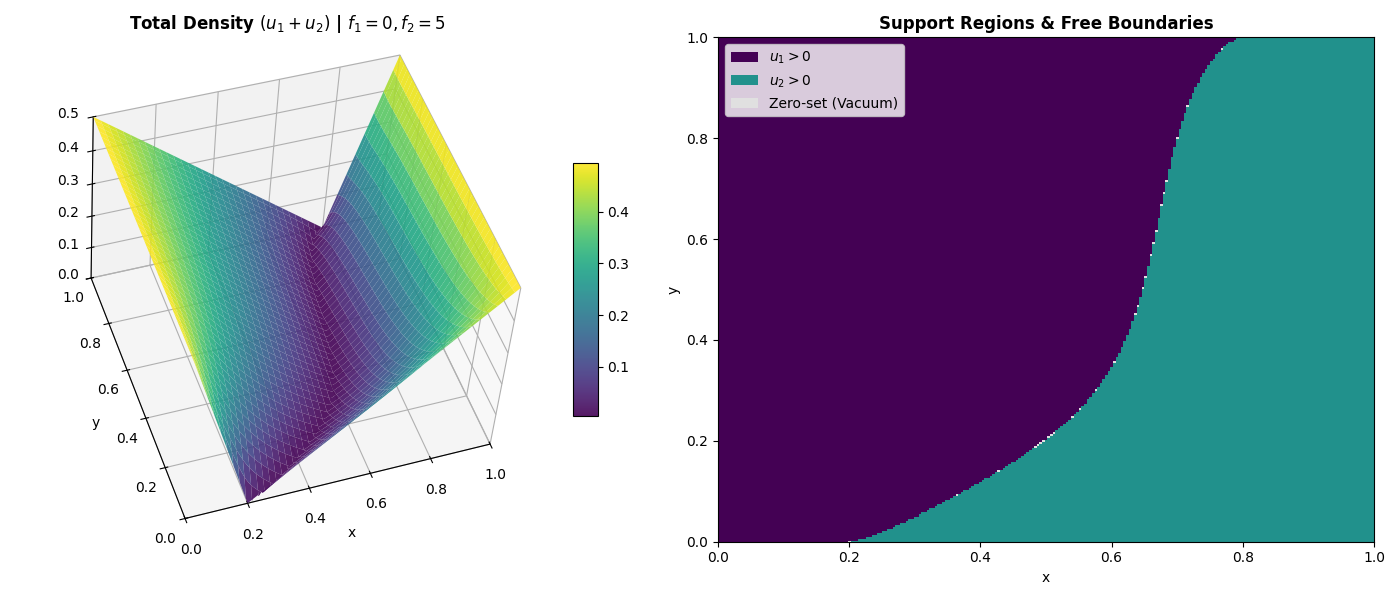}
\caption{Two-dimensional simulation with internal dynamics $f_1=0$ and $f_2=5$. The species meet along a common free boundary without the formation of a zero-set.}
  \label{2dfig}
\end{figure}

Figure \ref{2dfig} illustrates a scenario where the internal dynamics are relatively small ($f_1=0, f_2=5$). Here, we observe that no zero-set (vacuum region) forms between the species; the competing densities $u_1$ and $u_2$ meet directly along a common, continuous free boundary. 

\begin{figure}[!htbp]
    \centering
    \includegraphics[width=\textwidth]{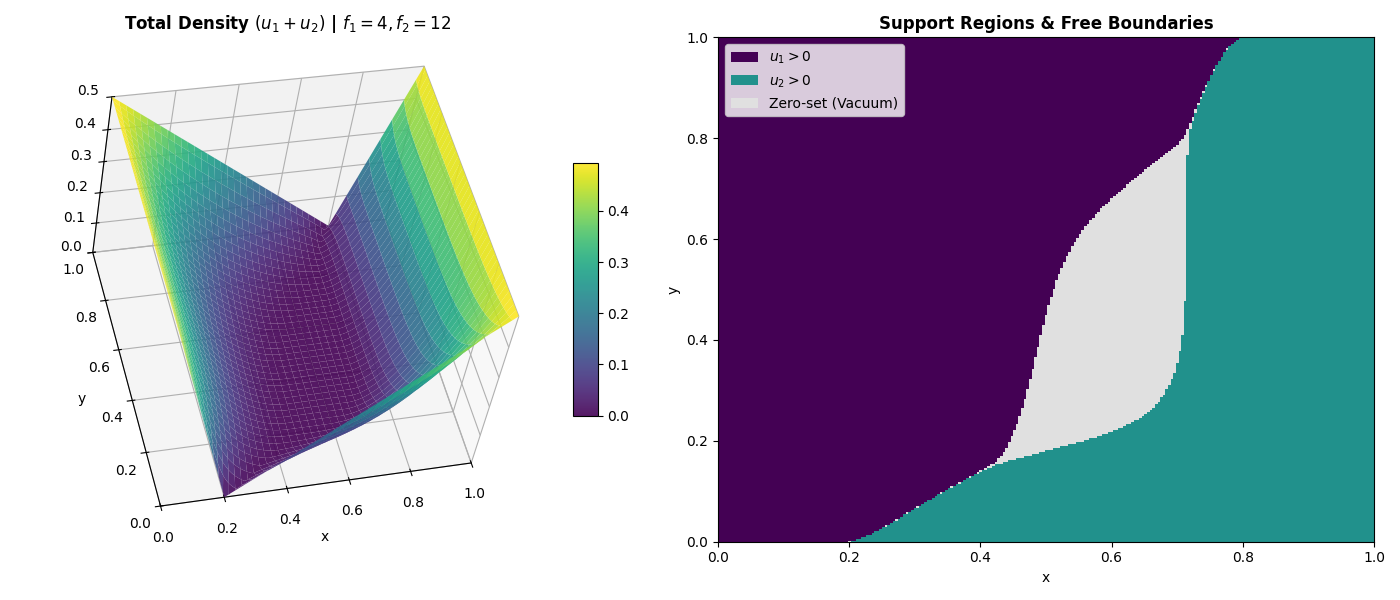}
 \caption{Two-dimensional simulation with internal dynamics $f_1=4$ and $f_2=12$. The high reaction penalties force a strict separation, creating a zero-set (vacuum region) between the two populations.}
  \label{2dfig2}
\end{figure}

In contrast, Figure \ref{2dfig2} demonstrates the effect of strongly penalizing internal dynamics ($f_1=4, f_2=12$). The high environmental cost forces both species to retreat from the competitive interface, resulting in the emergence of a distinct zero-set. This vacuum region physically separates the two populations, verifying that sufficiently large reaction terms strictly disconnect the active supports of the densities.

\bibliographystyle{acm}
\bibliography{conv}
\end{document}